\documentclass[pdflatex,iicol,sn-mathphys-num,10pt]{sn-jnl}
\usepackage{natbib}
\usepackage{graphicx}%
\usepackage{multirow}%
\usepackage{amsmath,amssymb,amsfonts}%
\usepackage{amsthm}%
\usepackage{mathrsfs}%
\usepackage[title]{appendix}%
\usepackage{xcolor}%
\usepackage{textcomp}%
\usepackage{manyfoot}%
\usepackage{booktabs}%
\usepackage{algorithm}%
\usepackage{algorithmicx}%
\usepackage{algpseudocode}%
\usepackage{listings}%
\usepackage{lscape}

\usepackage{xurl}
\usepackage{hyperref}

\theoremstyle{thmstyleone}%
\theoremstyle{thmstyletwo}%

\theoremstyle{thmstylethree}%

\begin{document}

\title[Article Title]{\centering{Residual Pseudospectra Reveal a Physics-Informed Koopman Backbone for\\Tropical Pacific Variability and ENSO Prediction}}


\author[1]{\fnm{Paula} \sur{Lorenzo-Sanchez}}\email{paula.lorenzo@cmcc.it}

\author*[2]{\fnm{Matthew} J. \sur{Colbrook}}\email{mjc249@cam.ac.uk}

\author[1,3]{\fnm{Antonio} \sur{Navarra}}\email{antonio.navarra@cmcc.it}

\affil[1]{\orgname{Centro Euromediterraneo sui Cambiamenti Climatic}, \orgaddress{\city{Bologna}}}

\affil[2]{\orgdiv{Department of Applied Mathematics and Theoretical Physics}, \orgname{University of Cambridge}}

\affil[3]{\orgdiv{Dipartimento di Scienze Geologiche, Biologiche e Ambientali}, \orgname{Universita’ di Bologna}}


\abstract{Tropical Pacific sea-surface-temperature (SST) variability spans interacting timescales, with the ENSO as its dominant interannual expression. Yet the dynamical structure organizing this variability and underpinning extended-range predictability remains difficult to extract from high-dimensional observations. Koopman operator learning offers spectral coordinates for nonlinear dynamics, yet finite geophysical records often produce dense, sampling-sensitive spectra whose physical content is ambiguous. We show that this apparent redundancy reflects coherent operator-level structure. Combining kernel Extended Dynamic Mode Decomposition with residual minimization and pseudospectral analysis, we use the Koopman eigenvalue relation as a physics-informed consistency test to organize learned spectra. Applied to ERA5 and HadISST tropical Pacific SST anomalies, the residual landscape identifies 19 robust residual-minimum frequencies with coherent spatial modes that persist across products and sampling realizations. Together, these modes define a compact Koopman backbone spanning low-frequency modulation through quasi-biennial components, including ENSO-band variability. The surrounding spectral cloud is structured by integer powers and nonlinear combinations of this backbone, forming a residual-ordered Koopman hierarchy. The backbone reconstructs substantial Niño3.4 variance and enables skillful out-of-sample forecasts, with greatest gains at 8–18-month leads. By embedding dynamical consistency into physics-informed operator learning, the framework turns opaque spectra into robust, interpretable and predictive representations of tropical Pacific variability.}




\maketitle


\begin{figure*}
 \centering
 \includegraphics[width=\textwidth]{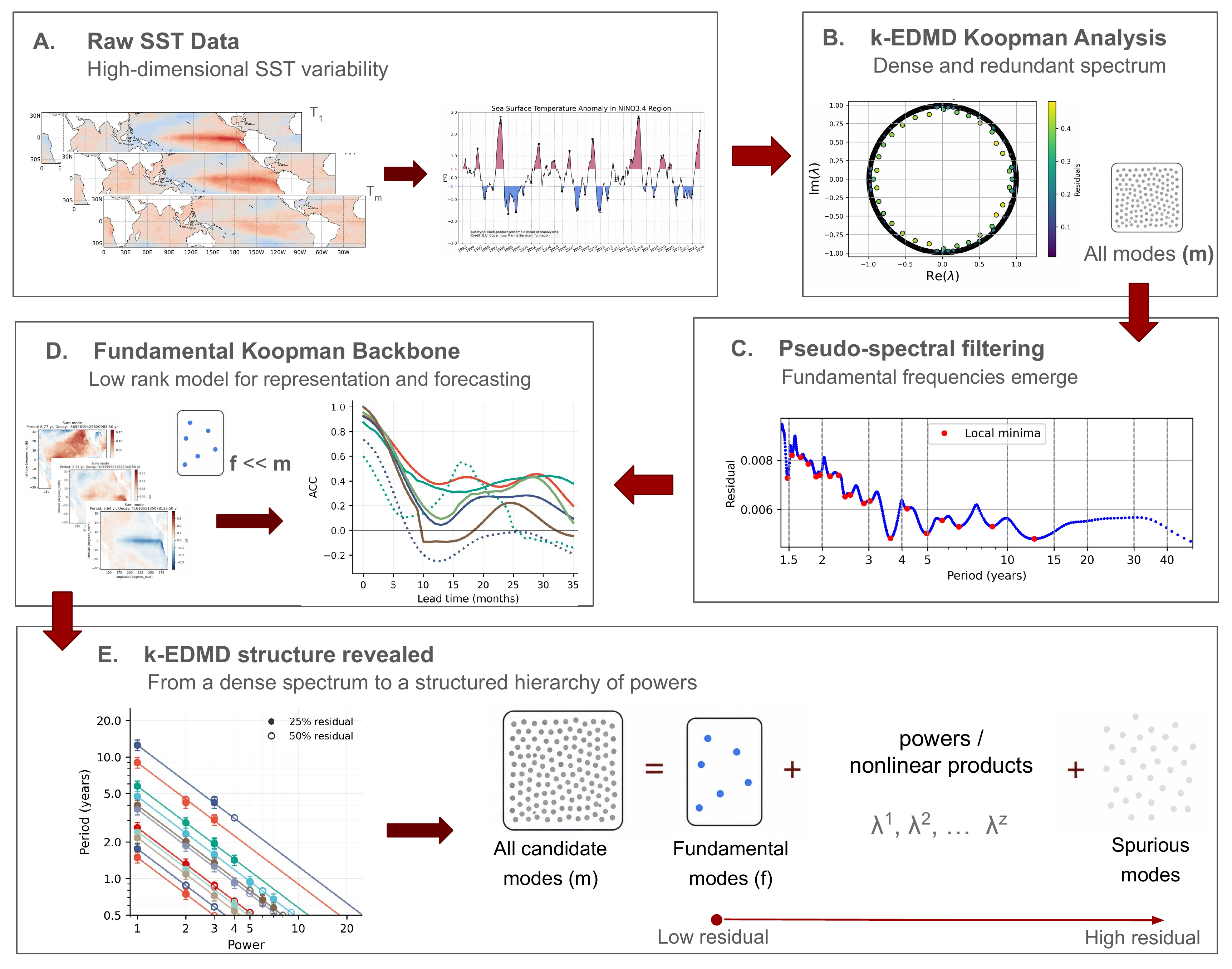}
 \caption{\textbf{Conceptual framework for extracting a Koopman backbone of ENSO variability.}
\textbf{A}, Monthly tropical Pacific sea-surface temperature fields are treated as high-dimensional observations of nonlinear climate variability associated with ENSO. 
\textbf{B}, Applying k-EDMD yields a dense cloud of Koopman eigenvalues near the unit circle, with substantial spectral redundancy and no immediately apparent low-dimensional organization. 
\textbf{C}, Pseudospectral analysis of the learned Koopman approximation along the unit circle identifies robust local minima in the residual norm. These minima define a compact set of well-resolved fundamental frequencies that persist across data products and sampling realizations, and are associated with coherent large-scale SST patterns. 
\textbf{D}, The corresponding fundamental modes form an interpretable Koopman backbone that reconstructs and forecasts ENSO variability. 
\textbf{E}, Analysis of the broader k-EDMD spectrum reveals a hierarchy in which many additional eigenpairs are organized as integer powers and nonlinear combinations of the fundamental components, while high-residual components are increasingly weakly resolved. 
Together, these steps transform a dense and redundant Koopman spectrum into a compact, physically interpretable, and predictive representation of tropical Pacific variability.} 
\label{Fig1}
\end{figure*}

Tropical Pacific sea-surface-temperature (SST) anomalies are a major source of global climate variability and predictability, coupling biennial, interannual, and decadal fluctuations across the ocean--atmosphere system \citep{Capotondi2023TPDV,Slawinska2017}. Its dominant interannual expression, the El Niño--Southern Oscillation (ENSO), reorganizes atmospheric circulation and shapes extreme weather, ecosystems, agriculture, water resources, and energy systems worldwide. The 2023--2024 El Niño highlighted this global reach, with impacts linked to drought and food insecurity in Southern Africa, as well as widespread marine heat stress and coral bleaching \citep{WMO2024ElNino, WFP2025SouthernAfrica, NOAA2024Bleaching}. Yet the dynamical structures that organize tropical Pacific variability, mediate nonlinear interactions across timescales, and constrain ENSO predictability remain incompletely resolved. This reflects a broader challenge for physics-informed machine learning in climate science: learning representations from high-dimensional data that are physically interpretable and dynamically consistent \citep{Heureux_2020, zhang_2022,GhilLucarini2020}.

As observational records, reanalysis products, and climate-model archives grow in length, resolution, and complexity, the bottleneck is shifting from data availability to dynamical interpretability. These datasets enable powerful data-driven discovery, but purely statistical decompositions and unconstrained latent representations often yield structures whose relation to physical mechanisms, predictive memory, and nonlinear climate dynamics is unclear. For climate applications, a useful learned representation must therefore do more than compress variance or optimize forecast skill: it must identify structures that are robust under sampling, interpretable in physical space, and consistent with the underlying dynamics \citep{LucariniChekroun2023, Lohmann2024Multistability}.

Koopman operator theory provides a natural bridge between data-driven learning and nonlinear dynamical systems. Rather than advancing state variables directly, it studies observables---functions or measurements of the state---that evolve under a linear, generally infinite-dimensional operator. Introduced by Koopman and von Neumann in the 1930s \citep{koopman_1931,koopman_1932}, this perspective has been transformed by modern data-driven approximation methods \citep{mezic_2005, Giannakis2019, Budisic2012AppliedKoopmanism, Brunton2022ModernKoopman, Colbrook2024MultiverseDMD}, enabling spectral analysis of complex systems across control, climate dynamics, epidemiology, neuroscience, non-autonomous systems, and interpretable machine learning \citep{Haggerty2023SoftRobots, Froyland2021, Orvieto2023RNN, Proctor2015DiseaseDMD, brunton2016, Lusch2018DeepKoopman, brunton2020machine, froyland2014computational, hogg2020, FroylandLloydSantitissadeekorn2010}. Its appeal for physics-informed machine learning is that the learned coordinates are not arbitrary latent variables: Koopman eigenfunctions and modes are tied to explicit temporal evolution. This makes Koopman analysis a principled route for extracting predictive and physically interpretable structure from nonlinear, partially observed climate data.

The Koopman spectrum provides a temporally ordered representation of nonlinear dynamics with direct dynamical interpretation \citep{mezic_2005,Rowley2009,Mezic2013,Williams2015Kernel}. Koopman eigenvalues specify characteristic timescales, including oscillation frequencies and decay rates; eigenfunctions define intrinsic coordinates with simple temporal evolution; and Koopman modes project these temporal signatures onto spatial patterns in the observed physical variables. This differs fundamentally from covariance- or variance-based decompositions such as EOF analysis, which rank structures by instantaneous amplitude. For climate data, this distinction is critical: dynamically persistent or predictive structures need not be the most energetic, and high-variance patterns need not encode long-term memory. Koopman analysis therefore offers a physics-informed framework for separating interacting timescales, identifying structures with dynamical persistence, and constructing reduced-order models for reconstruction and prediction \citep{navarra_2021, Froyland2024NonAutonomous}.

Data-driven Koopman approximations, including Dynamic Mode Decomposition (DMD), Extended DMD (EDMD), and kernel EDMD (k-EDMD), have made operator-theoretic learning practical for high-dimensional empirical systems \citep{Klus_2016,Klus2018}. Yet in long geophysical records, this practicality exposes a new challenge: the learned spectrum is often dense, redundant, and sensitive to finite-sample effects. This issue is amplified in snapshot-based k-EDMD, where longer training records improve the sampling of slow climate variability but also increase the dimension of the spectral problem. In applications to tropical Pacific SST, this can lead to hundreds of eigenpairs clustered near the unit circle, a structure consistent with persistent oscillatory dynamics but difficult to interpret mode by mode \citep{Lorenzo2025,Navarra2024}. Small changes in data product, sampling window, or numerical approximation can reorganize the discrete eigenvalue cloud, obscuring which components represent robust climate dynamics and which reflect spectral redundancy or weak finite-data resolution \citep{brunton2016,brunton2017,Williams2015,lewin2010,mezic_2005,Lorenzo2025}.

The central problem is therefore not merely the computation of Koopman spectra, but their interpretation. Existing selection strategies typically rank individual modes by growth rate, frequency, variance contribution, or predictive skill \citep{Tu2014,pan2019efficient,navarra_2021}. These criteria are useful, but they treat eigenpairs as isolated candidates and do not resolve the organization of the learned spectrum as a whole. Dense finite-data spectra may contain robust oscillatory components, finite-data approximations of mixed or continuous spectral content, algebraic products of Koopman eigenfunctions, redundant representations, and weakly resolved numerical structures \citep{GiannakisValva2024}. Distinguishing among these possibilities is essential for physics-informed prediction, particularly in the tropical Pacific, where ENSO emerges from nonlinear interactions across seasonal, interannual, and decadal timescales.

Earlier reduced-order models of tropical Pacific variability established that compact dynamical representations can be highly predictive. Linear Inverse Models (LIMs), for example, approximate climate anomalies as a stochastically forced linear system whose eigenmodes capture dominant patterns of variability, growth, and predictability \citep{Penland1995,Newman2011}. These models showed that a small set of linear modes can reproduce substantial ENSO variability and provide a strong benchmark for seasonal prediction. However, their prescribed linear dynamics limit their ability to represent nonlinear interactions and state-dependent dynamics that shape tropical Pacific variability.

Recent nonlinear operator-theoretic and related spectral methods have shown that climate variability contains coherent multiscale structure. Using nonlinear Laplacian spectral analysis, Slawinska and Giannakis \citep{Slawinska2017} recovered a hierarchy of Indo-Pacific SST modes, including ENSO, ENSO--annual-cycle combinations, tropospheric biennial oscillation-related variability, and decadal Pacific structures. Froyland et al. \citep{Froyland2021} subsequently showed that Koopman and transfer-operator eigenfunctions provide dynamically consistent coordinates for the ENSO lifecycle, identifying a fundamental ENSO mode, a 3-year ENSO component, and combination modes between ENSO and the annual cycle. Related operator-theoretic response approaches have used Kolmogorov modes to connect natural variability and forced response in stochastic ENSO-type models \citep{ChekrounZagliLucarini2025}.  Koopman and related kernel methods have also supported Niño3.4 prediction from tropical Pacific SST \citep{Berry_2015,wang2020,navarra_2021,Lorenzo2025}, isolated periodic atmospheric modes and nonlinear interactions in stratospheric winds \citep{Valva2025QBO}, and enabled long-term prediction of sea-ice variability from satellite observations \citep{hogg2020,ColbrookMezicStepanenko2026}. These studies establish operator-theoretic learning as a powerful framework for discovering climate-relevant modes, nonlinear spectral interactions, and predictive coordinates. What remains unresolved is how to interpret the dense finite-data spectra produced by nonlinear Koopman approximations themselves: whether they contain an intrinsic organization that separates robust primary dynamics from Koopman-algebraic relationships, redundancy, and weakly resolved spectral content.

Residual and pseudospectral Koopman methods provide the operator-level diagnostics needed for this classification. Residual Dynamic Mode Decomposition (ResDMD) quantifies how closely computed eigenpairs satisfy the Koopman eigenvalue problem, helping distinguish well-resolved spectral information from numerical artefacts or weakly supported modes \citep{ColbrookTownsend2024,colbrook2022}. One may also study the SVD truncation often involved in DMD \citep{drmac2018,drmac2023}. Pseudospectral analysis adds a complementary perspective by measuring spectral sensitivity to perturbations and identifying regions of the learned operator that are robustly resolved \citep{TrefethenEmbree2005,ColbrookTownsend2024,GiannakisValva2024}. Together, these tools make it possible to move beyond heuristic mode-by-mode selection toward a spectrum-level, physics-informed assessment of Koopman reliability, redundancy, and dynamical interpretability.

Here we use residual and pseudospectral structure as an operator-level organizing principle for dense k-EDMD spectra. Rather than selecting isolated ENSO-related modes or focusing on prescribed spectral interactions, we make the finite-data Koopman spectrum itself the object of analysis. Using detrended monthly SST anomalies with the climatological annual cycle removed, we focus on quasi-biennial to decadal tropical Pacific anomaly dynamics. Residual minima identify robust candidate frequencies, kernel-independent spatial reconstructions test their physical coherence, and the remaining eigenpairs are classified according to whether they are consistent with harmonics, nonlinear combinations, redundancies, or weakly resolved finite-data structure. This converts Koopman learning from heuristic mode selection into a physics-informed spectral-organization problem, linking numerical robustness, physical interpretability, nonlinear spectral relationships, and predictive relevance within a single hierarchy.

Figure~\ref{Fig1} summarizes the workflow. Starting from high-dimensional tropical Pacific SST data, snapshot k-EDMD produces a dense near-unit-circle spectrum whose physical organization is not apparent from the eigenvalue cloud alone. Residual minimization and pseudospectral analysis identify well-resolved frequencies and their associated reconstructed SST patterns define a compact Koopman backbone of tropical Pacific variability, within which ENSO emerges as the dominant interannual component. The remaining eigenpairs are then interpreted relative to this backbone: many align with harmonics of the robust components, whereas high-residual eigenpairs form increasingly weakly resolved finite-data structure. The result is a spectrum-level decomposition of the learned Koopman operator into robust primary components, dynamically related secondary structure, and weakly supported residual modes.

Using tropical Pacific SST anomalies from ERA5 and HadISST, we validate this spectrum-level interpretation through four complementary tests. First, the residual-minimum frequencies and reconstructed spatial modes are reproducible across data products and sampling realizations, and align with physically interpretable tropical Pacific and ENSO-related SST anomaly patterns. Second, the resulting 19-component Koopman backbone reconstructs approximately 92\% of training-period Niño3.4 variance while retaining only a small, interpretable set of modes. Third, reduced-order forecasts based on this backbone retain out-of-sample Niño3.4 skill, with the strongest advantages at 8--18-month lead times, indicating that the selected modes encode dynamical memory rather than only in-sample reconstruction skill. Fourth, many remaining eigenpairs are organized as integer powers or nonlinear combinations of the backbone, revealing a Koopman-algebraic hierarchy within the dense k-EDMD spectrum. Together, these tests show that residual and pseudospectral analysis can convert spectral redundancy from an interpretability bottleneck into a physics-informed representation of nonlinear climate variability that is robust, interpretable, and predictive.
\subsection*{Results}

\subsubsection*{Pseudospectral residuals reveal stable structure hidden within dense Koopman spectra}

\begin{figure*}
 \centering
 \includegraphics[width=\textwidth]{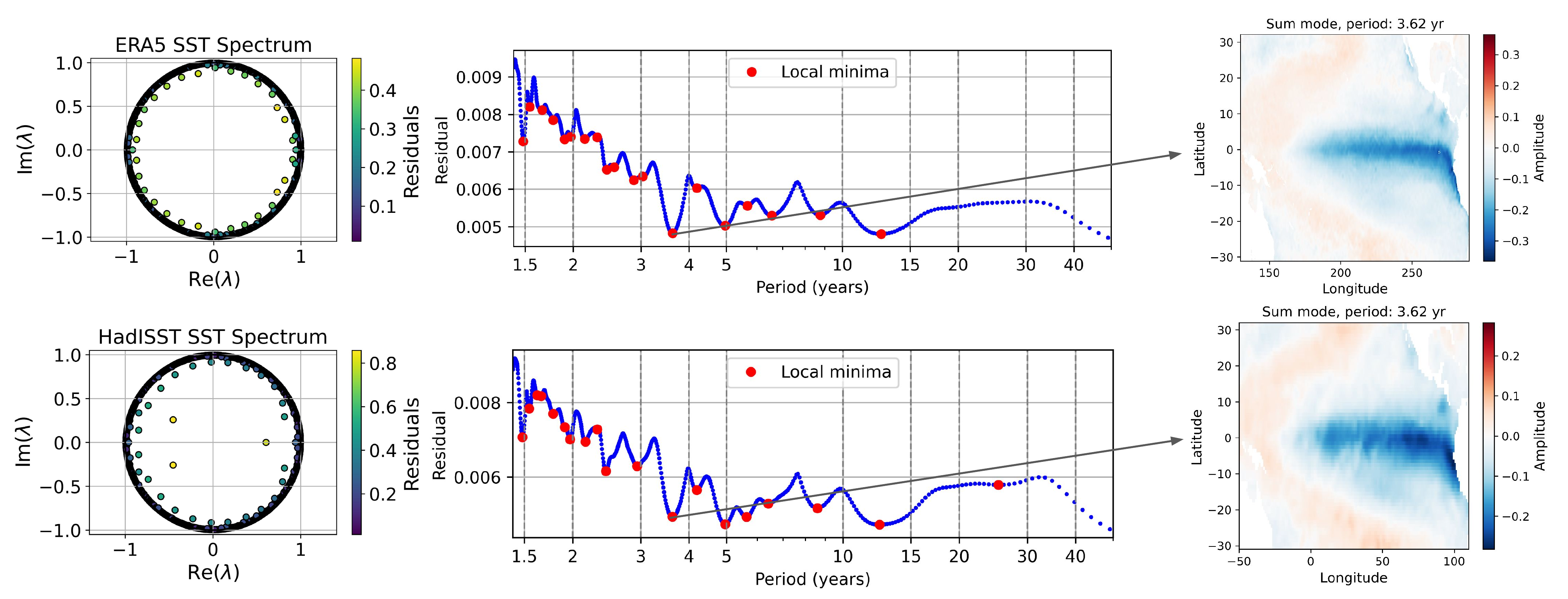}
 \caption{\textbf{k-EDMD spectra and pseudospectral residual structure for tropical Pacific SST.}
\textbf{Left}, k-EDMD eigenvalue spectra computed from monthly tropical Pacific SST anomalies in ERA5 (top) and HadISST (bottom) over the 1940--2010 training period. Colours denote the residual norm associated with each computed eigenpair.
\textbf{Middle}, pseudospectral residual evaluated along the unit circle as a function of positive frequency. Local residual minima (red markers) identify periods at which the Koopman approximation is robustly resolved; we refer to these as \textit{fundamental periods}.
\textbf{Right}, real part of the spatial Koopman modes reconstructed at a dominant pseudospectral minimum with period $\approx 3.6$ years using the kernel-independent summation procedure. The reconstructed patterns display canonical ENSO-like SST structure and are highly consistent between ERA5 and HadISST, despite differences in spatial resolution.} 
\label{Fig2}
\end{figure*}

Applying k-EDMD to monthly tropical Pacific SST anomaly fields over the 1940--2010 training period yields dense Koopman spectra for both ERA5 and HadISST (Fig.~\ref{Fig2}, left column). Each data product produces approximately 800 eigenpairs, with most eigenvalues lying close to the unit circle. This pattern is consistent with weakly damped, persistent tropical Pacific variability, including interannual ENSO-related oscillations \citep{Navarra2024,Lorenzo2025}. However, this spectral density creates an interpretability bottleneck: many eigenpairs share similar frequencies and decay rates, obscuring which components correspond to distinct physical dynamics and which reflect redundancy or weak spectral resolution.

The overall spectral envelopes are similar in ERA5 and HadISST, indicating that both products capture comparable large-scale tropical Pacific variability. At the level of individual eigenpairs, however, the spectra differ substantially. Eigenvalues shift in phase and magnitude, residual norms vary, and one-to-one matching across products is not well defined, even among low-residual eigenpairs near the unit circle. Residual norms therefore provide valuable information about spectral reliability, but residual ordering of the discrete eigenpairs alone does not give a unique or robust criterion for selecting physically meaningful modes.

The same ambiguity appears under small perturbations of the training window. For 65-year ERA5 windows shifted by only one year, the detailed eigenvalue distribution and the dominant peaks of the eigenvalue-density estimate change substantially (Supplementary Fig.~S1). Raw spectral-density peaks therefore do not, by themselves, define stable dynamical frequencies. By contrast, pseudospectral residual curves computed over the same windows remain highly consistent, with local minima occurring at nearly identical periods. Thus, the instability lies mainly in the placement and density of the discrete k-EDMD eigenpairs, whereas the operator-residual landscape reveals a more stable spectral organization.

Figure~\ref{Fig2} summarizes this contrast. Although the raw eigenvalue clouds are dense and product-dependent, the pseudospectral residual along the unit circle contains reproducible local minima in both ERA5 and HadISST. These minima identify frequencies at which the learned Koopman approximation is robustly supported, providing the basis for the pseudospectral filtering step below.

\subsubsection*{Pseudospectral filtering identifies robust fundamental frequencies}

To extract robust dynamical structure from the dense spectra, we evaluate the pseudospectral residual of the learned Koopman approximation along the unit circle (Fig.~\ref{Fig2}, central column). In both ERA5 and HadISST, the residual landscapes contain well-defined local minima at nearly identical frequencies, despite substantial differences in the corresponding discrete eigenvalue clouds. Across data products and training windows, we identify 19 such residual-minimum components within the approximately 1.5-year-to-decadal period range. We refer to these as the fundamental frequencies of the finite-data tropical Pacific Koopman representation.

These frequencies span the principal timescales of tropical Pacific SST variability, from low-frequency modulation through the classical ENSO band to shorter quasi-biennial components. Several residual minima lie within the ENSO band, including a prominent component near 3.6 years. Reconstructing the associated spatial Koopman modes using the kernel-independent summation procedure in Eq.~\ref{summation_modes} yields coherent SST structures that are reproducible across ERA5 and HadISST, despite differences in spatial resolution and data construction. The 3.6-year mode (Fig.~\ref{Fig2}, right column) displays a canonical ENSO-like pattern, with anomalies concentrated in the central and eastern equatorial Pacific and weaker westward extensions.

The full set of reconstructed modes falls into a small number of physically interpretable families (Supplementary Table~1; Supplementary Figs.~S2 and~S3). The longest-period components show broad basin-scale structures with substantial off-equatorial amplitude and meridional asymmetry, consistent with low-frequency tropical Pacific modulation. ENSO-band components, with periods of approximately 3--6 years, are concentrated along the equatorial cold tongue. Shorter-period components retain an equatorial Pacific signature but become increasingly localized, consistent with faster ENSO-related and quasi-biennial variability. This multiscale organization is broadly consistent with earlier operator-theoretic analyses of Indo-Pacific and ENSO variability, which identified ENSO-band, quasi-biennial and decadal Pacific components \citep{Slawinska2017,Froyland2021}.

\begin{figure*}
 \centering
 \includegraphics[width=0.9\textwidth]{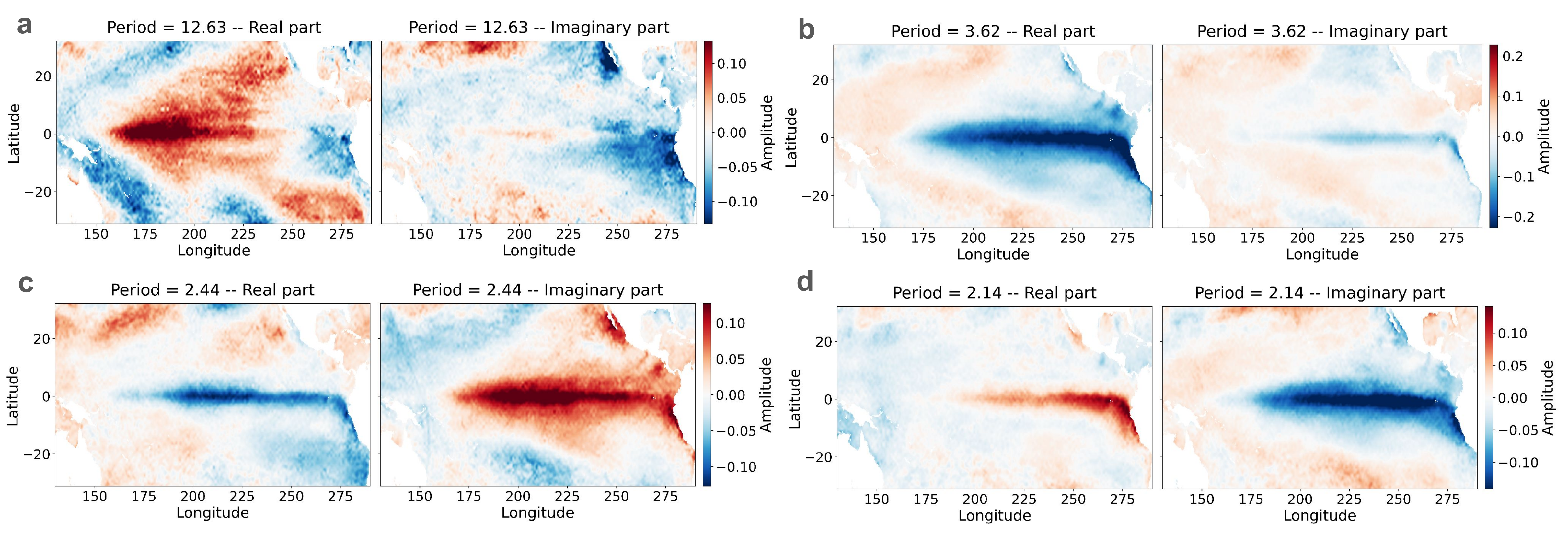}
 \caption{\textbf{Representative fundamental Koopman modes of tropical Pacific SST variability.}
Real and imaginary parts of four ERA5 complex-valued modes reconstructed at pseudospectral residual minima. For each mode, the two parts correspond to quadrature phases of a single oscillatory pattern; the absolute sign and phase are arbitrary.
\textbf{a}, The $\sim 12.6$-year mode shows broad basin-scale SST structure with substantial off-equatorial amplitude, consistent with low-frequency tropical Pacific modulation.
\textbf{b}, The $\sim 3.6$-year ENSO-band mode is dominated by central and eastern equatorial Pacific cold-tongue anomalies.
\textbf{c}, The $\sim 2.4$-year mode has a central-Pacific ENSO-like structure, with anomalies concentrated in the central equatorial Pacific.
\textbf{d}, The $\sim 2.1$-year mode shows a more localized equatorial pattern, consistent with quasi-biennial ENSO-related variability.
The full set of 19 fundamental modes is shown in Supplementary Figs.~S2 and~S3.} 
\label{Fig3}
\end{figure*}

\begin{figure*}
 \centering
 \includegraphics[width=\textwidth]{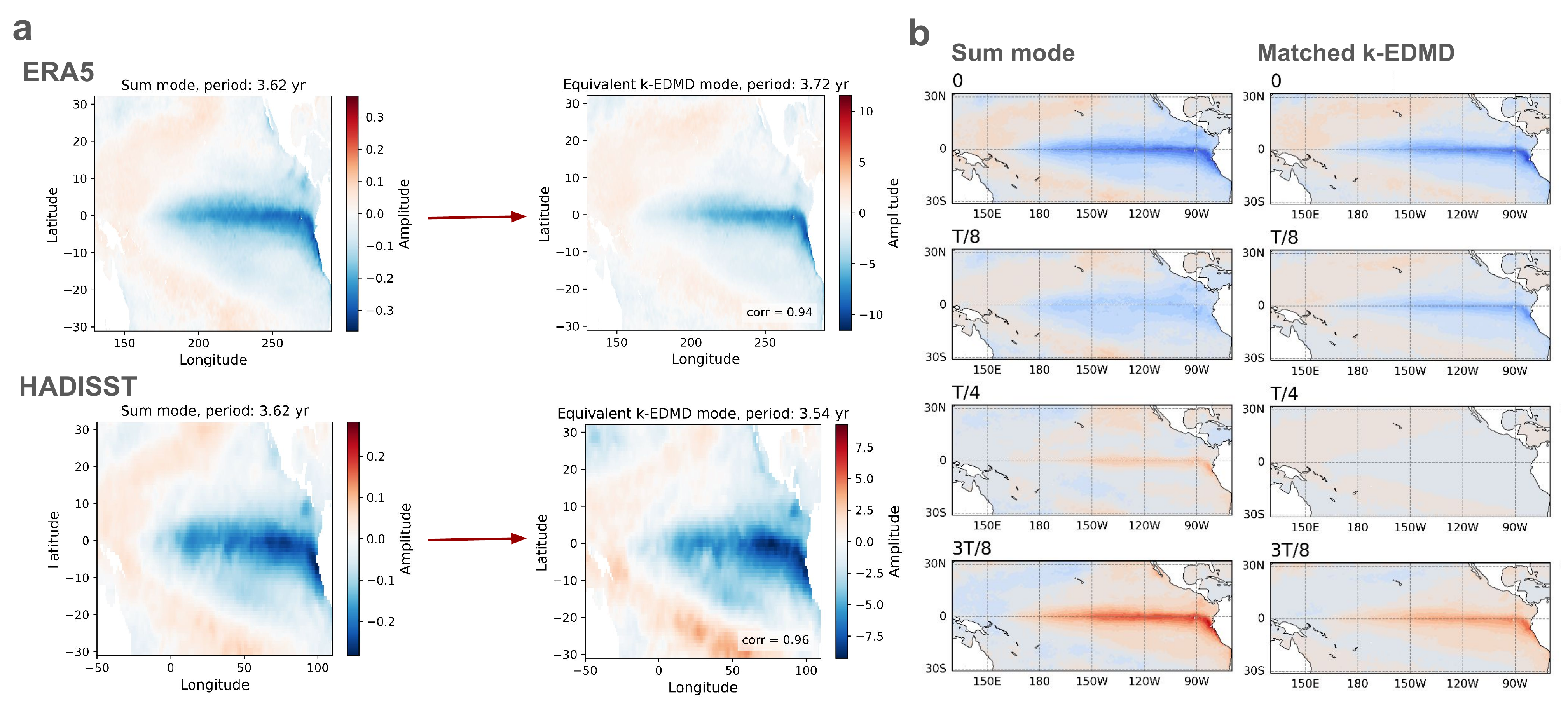}
 \caption{
 \textbf{Pseudospectral fundamental modes and their k-EDMD spectral counterparts.}
\textbf{a}, Pseudospectral fundamental mode and matched k-EDMD eigenmode for the dominant ENSO-like period of $\sim 3.6$ years. Results are shown for ERA5 (top) and HadISST (bottom). Left maps show the real part of the spatial mode reconstructed from the pseudospectral residual minimum; right maps show the real part of the best-correlated k-EDMD eigenmode identified in the original spectrum. In both data products, the k-EDMD counterpart reproduces the large-scale ENSO-like SST structure with high spatial correlation (lower right corner).
\textbf{b}, Phase evolution of the same pseudospectral mode and its matched k-EDMD counterpart for ERA5. Maps show the first half of the oscillatory cycle at quarter-cycle intervals; the second half is omitted because it repeats the sequence with reversed sign. The close agreement shows that, when a matched k-EDMD counterpart exists, it captures both the spatial structure and phase evolution of the corresponding ENSO-like oscillation.
 }
\label{Fig4}
\end{figure*}

Figure~\ref{Fig3} shows four representative modes from this hierarchy. Because the modes are complex-valued and occur in conjugate pairs, each pair represents an oscillatory SST pattern with two quadrature phases. For a mode with period \(T\), the real and imaginary components correspond to phases separated by \(T/4\), up to an arbitrary phase convention; after half a cycle, the same structures reappear with reversed sign. Physical interpretation therefore rests on the coherent spatial pattern and its phase evolution, rather than on the absolute sign or phase of an individual component. This evolution is shown explicitly for the \(\sim 3.6\)-year mode in Fig.~\ref{Fig4}b.

The \(\sim 12.6\)-year mode represents the low-frequency part of the Koopman backbone. Its spatial pattern spans much of the tropical Pacific, with substantial off-equatorial amplitude and meridional asymmetry. Similar broad structures appear among other longer-period modes in Supplementary Fig.~S2, particularly near 8.8, 6.6, and 5.7 years. Given the limited number of decadal cycles in the training record, we interpret this family cautiously as low-frequency modulation of tropical Pacific SST variability, potentially related to PDO/IPO-like modulation of ENSO and broadly comparable to low-frequency Pacific modes identified in previous operator-theoretic studies \citep{Slawinska2017,navarra_2021,Capotondi2023TPDV}. Its off-equatorial structure is also suggestive of subtropical--tropical coupling pathways \citep{ChiangVimont2004,Vimont2003}, although the SST-only patterns should not be interpreted as a direct Pacific Meridional Mode analogue.

The $\sim 3.6$-year mode is the clearest Koopman representation of canonical ENSO variability. Its SST anomalies are concentrated in the central and eastern equatorial Pacific, with largest amplitudes along the cold tongue. The quadrature phase captures the phase-shifted continuation of this pattern, showing that the mode encodes an evolving ENSO-like oscillation rather than a static anomaly field. Its period and cold-tongue structure place it in the same broad ENSO-band family as the fundamental ENSO components identified by Slawinska and Giannakis, while here the analysis is restricted to SST anomalies and seasonal phase-locking is not represented explicitly. In the full set of modes, the components with periods between roughly 3 and 6 years form a related ENSO-band family, indicating that the Koopman backbone does not reduce ENSO variability to a single frequency but captures a small set of robust interannual components.

At shorter periods, the \(\sim 2.4\)-year mode retains a pronounced equatorial signature but is more strongly weighted toward the central Pacific than the \(\sim 3.6\)-year mode. This spatial structure suggests a faster ENSO-related component with characteristics reminiscent of Central Pacific ENSO-like variability, in which SST anomalies are organized preferentially around the central equatorial Pacific rather than extending as broadly across the eastern cold-tongue region \citep{Capotondi2015ENSODiversity}. Similar central-equatorial structures appear among nearby modes in Supplementary Fig.~S2, including components with periods near \(\sim 3\) and \(\sim 2.4\) years. We therefore interpret it as part of the transition from canonical ENSO-band variability toward faster interannual components.

The $\sim 2.1$-year mode marks the transition toward higher-frequency tropical Pacific variability. Its pattern remains equatorially organized but is more localized, resembling other short-period components in the fundamental set, including modes near $\sim 1.90$ and $\sim 1.48$ years (Supplementary Figs. S2 and S3). Although its period lies in the quasi-biennial range previously identified in equatorial Pacific variability \citep{Jiang1995}, its spatial structure is more localized than the broader Indo-Pacific TBO-related patterns identified by Slawinska and Giannakis \citep{Slawinska2017}. We therefore interpret it cautiously as a short-period ENSO-related anomaly component, which may either represent quasi-biennial tropical Pacific variability or arise as a harmonic or nonlinear product of lower-frequency Koopman components, as examined in Section~\ref{sec:hierarchy}.

Supplementary Table~1 reports the 19 ERA5 fundamental periods, their HadISST counterparts, and their assigned dynamical families. We match a HadISST mode to an ERA5 mode when its period lies within 10\% of the ERA5 period and its spatial correlation exceeds 0.7, with correlations reported in parentheses. Most ERA5 fundamental modes have a HadISST counterpart, indicating that the selected frequencies and spatial structures are reproducible across data products. This reproducibility is strongest for the ENSO-band modes and for several longer-period equatorial components.

Although the fundamental modes are identified from pseudospectral residual minima, many have nearby counterparts in the original k-EDMD spectrum. We refer to these as matched modes: k-EDMD eigenpairs whose periods and spatial patterns closely correspond to a pseudospectral fundamental mode. These matched modes provide a link between the robust pseudospectral representation and the discrete eigenpairs returned by k-EDMD.

Figure~\ref{Fig4} illustrates this correspondence for the dominant ENSO-like component with period \(\sim 3.6\) years. In both ERA5 and HadISST, the pseudospectral reconstruction yields a coherent equatorial Pacific SST pattern, and a closely matching k-EDMD eigenmode is present in the original spectrum (Fig.~\ref{Fig4}a). Spatial correlations of 0.94 for ERA5 and 0.96 for HadISST show that the two representations recover the same large-scale ENSO-like structure. For ERA5, Fig.~\ref{Fig4}b further shows that this correspondence extends to the phase evolution of the oscillation.

Across the full fundamental set, most pseudospectral modes have matched k-EDMD counterparts with low residuals and high spatial correlations (Supplementary Table~1). These counterparts show that robust structures are embedded within the dense eigenvalue cloud. However, the discrete representatives are not fully stable: small training-window shifts can cause matched eigenpairs to move, merge or lose an identifiable counterpart. Thus, matched k-EDMD modes are useful spectral representatives when present, but the pseudospectral residual minima provide the more stable definition of the fundamental Koopman components.

\subsubsection*{A 19-mode Koopman backbone reconstructs most Niño3.4 variability}

\begin{figure*}
 \centering
 \includegraphics[width=0.9\textwidth]{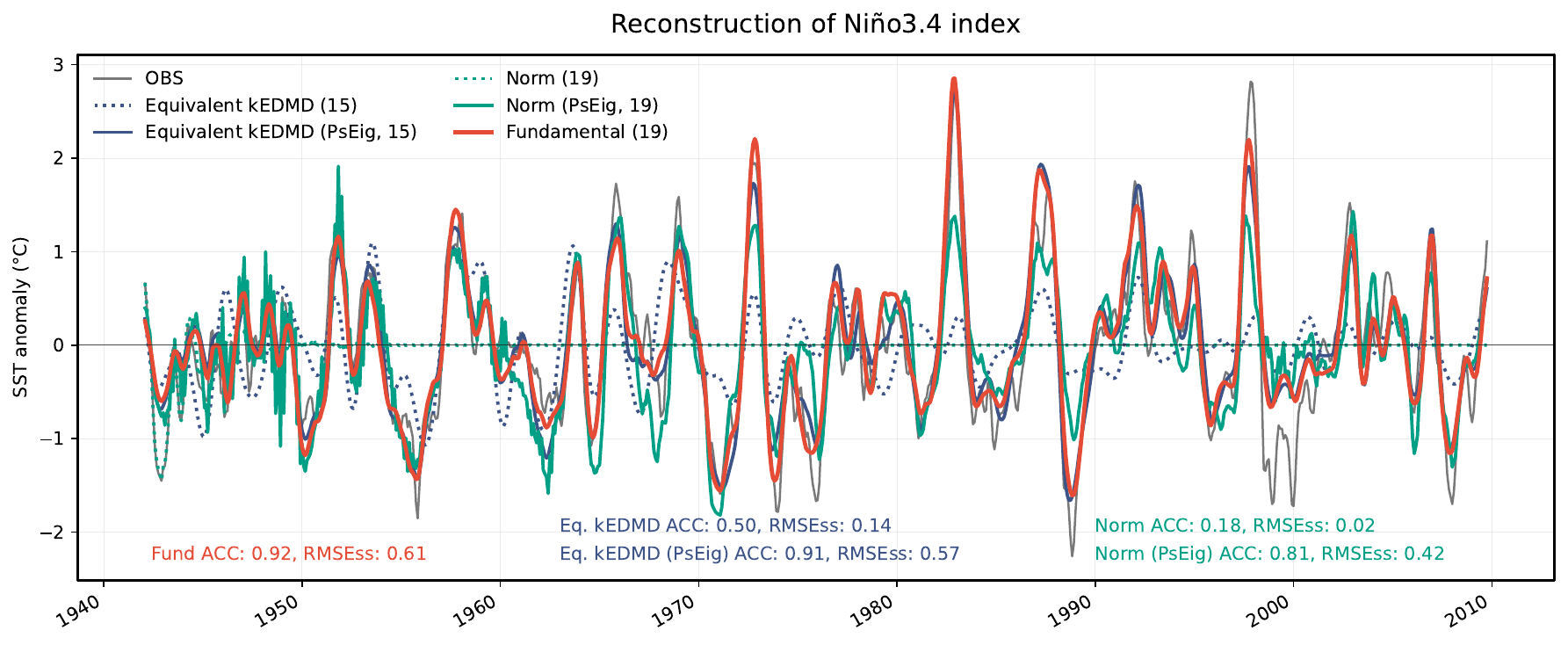}
 \caption{
 \textbf{Reduced-order reconstruction of the Niño3.4 index from different Koopman mode selections.}
The observed Niño3.4 SST anomaly index (grey) is compared with reconstructions obtained from different Koopman modal subsets. The reconstruction using the 19 pseudospectral fundamental modes (red) captures approximately 92\% of the observed variance. Reconstructions using the matched k-EDMD modes identified in the original spectral cloud explain substantially less variance when standard k-EDMD eigenfunctions are used (blue dotted; $\sim 50\%$), but improve markedly when the same modes are represented with their associated pseudoeigenfunctions (blue solid). As an amplitude-based benchmark, reconstructions using the 19 modes with the largest singular-value norms (teal blue dotted and solid) perform worse than the pseudospectral fundamental modes, despite selecting energetically dominant components. Correlation and RMSE skill scores for each reconstruction are shown in the lower part of the panel.
 } 
\label{Fig5}
\end{figure*}

The pseudospectral fundamental modes define a compact set of robust frequencies with coherent tropical Pacific SST structure. To test whether this backbone captures dynamically relevant variability, we reconstruct the Niño3.4 index using reduced-order Koopman representations built from different mode subsets (Fig.~\ref{Fig5}).

Using the 19 pseudospectral fundamental modes, the reduced Koopman reconstruction captures approximately 92\% of the training-period Niño3.4 variance (red line, Fig.~\ref{Fig5}). Thus, although k-EDMD returns hundreds of eigenpairs, the dominant interannual ENSO signal is well represented by a much smaller set of robust and physically interpretable components. This supports the interpretation of the residual-minimum modes as a compact Koopman backbone of tropical Pacific SST variability.

We next compare this reconstruction with one based on the matched k-EDMD modes, whose periods and spatial patterns closely match the pseudospectral fundamental modes. When these modes are reconstructed using their standard k-EDMD eigenfunctions, they explain only about 50\% of the Niño3.4 variance (blue dotted line, Fig.~\ref{Fig5}). Replacing these eigenfunctions with the corresponding pseudoeigenfunctions yields a reconstruction comparable to that from the pseudospectral fundamental modes (solid blue line). The loss of reconstruction skill is therefore not caused primarily by incorrect frequencies or spatial patterns, but by instability in the standard k-EDMD eigenfunction representation.

To test whether reconstruction skill is driven simply by variance retention, we also reconstruct Niño3.4 using the 19 k-EDMD modes with the largest singular-value norms, a standard heuristic in reduced-order modelling. These energetic modes explain only 18\% of the variance when reconstructed with standard k-EDMD eigenfunctions and 81\% when reconstructed with the corresponding pseudoeigenfunctions (teal lines, Fig.~\ref{Fig5}). Thus, large mode norm alone does not identify the most dynamically informative components: the pseudospectral fundamental modes recover substantially more Niño3.4 variability despite being selected by residual robustness rather than amplitude. Reconstructions from 19 randomly selected modes perform worse, explaining on average only 2\% of the variance with standard k-EDMD eigenfunctions and 64\% with pseudoeigenfunctions over 100 random subsets.

Together, these comparisons show that pseudospectral filtering identifies a compact set of dynamically informative Koopman modes. The fundamental modes are not merely a visually coherent subset of the spectrum; they capture the dominant Niño3.4 variability while suppressing much of the redundancy and eigenfunction instability present in the full k-EDMD spectral cloud.

\subsubsection*{The remaining k-EDMD spectrum forms a Koopman-algebraic hierarchy}
\label{sec:hierarchy}

The 19-mode backbone captures much of the dynamically relevant tropical Pacific variability, but k-EDMD returns hundreds of additional eigenpairs. This raises a central interpretability question: do the remaining eigenpairs represent independent dynamical degrees of freedom, or are many of them generated by the algebraic structure of Koopman eigenfunctions?

\begin{figure*}
 \centering
 \includegraphics[width=\textwidth]{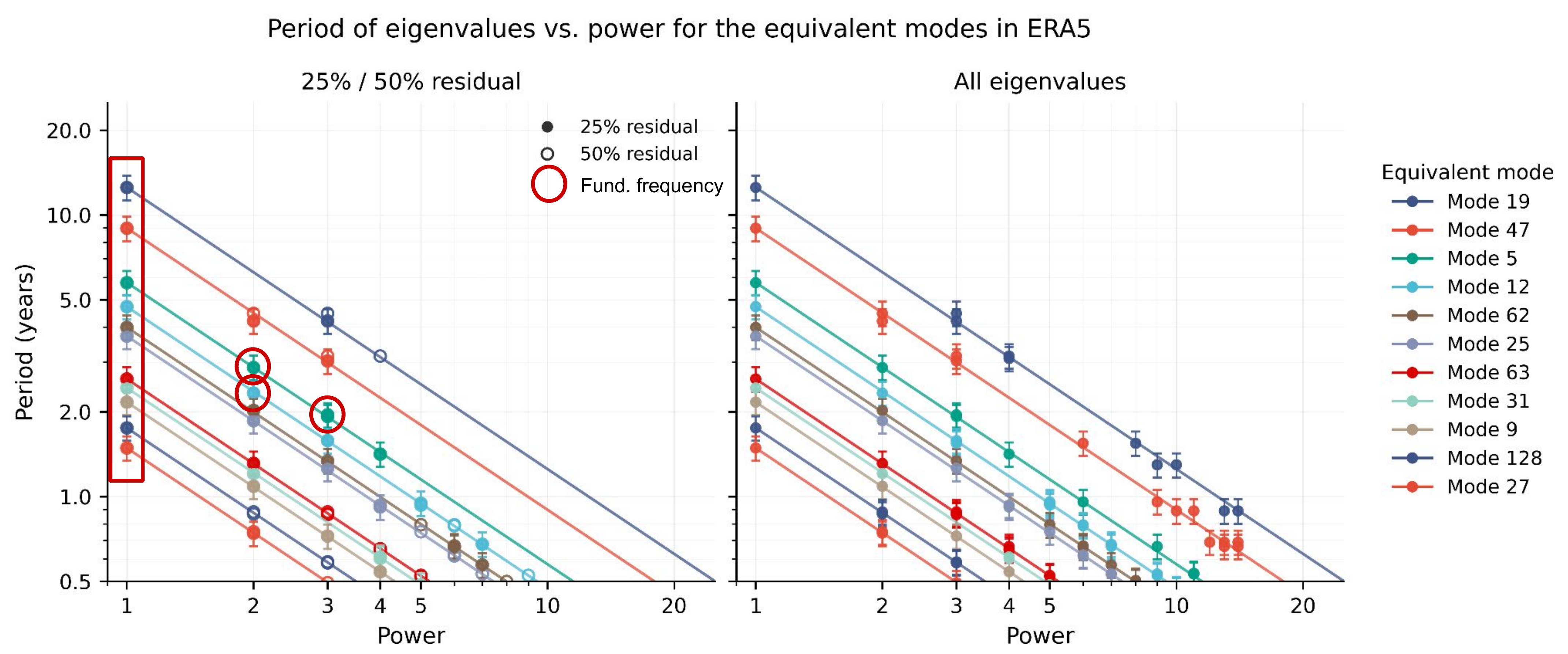}
 \caption{
 \textbf{Hierarchical organization of k-EDMD eigenpairs generated by fundamental modes (ERA5).}
Eigenvalue-implied periods of k-EDMD eigenpairs matched to integer powers of the equivalent fundamental modes, plotted as a function of power order. Each curve corresponds to one matched k-EDMD seed mode associated with a pseudospectral fundamental frequency. Labels indicate the eigenpair index after ordering by residual norm. Markers denote eigenpairs whose eigenvalue-implied periods and eigenfunction subspaces match the corresponding power within the prescribed tolerances. \textbf{Left}, hierarchy restricted to the lowest-residual 50\% of the k-EDMD spectrum. Filled markers denote eigenpairs in the lowest-residual 25\%, while open markers denote eigenpairs between the 25th and 50th residual percentiles. Error bars show the \(\pm 10\%\) period tolerance used for power matching. Red outlines highlight the matched fundamental modes.
\textbf{Right}, same analysis using the full k-EDMD spectrum. Including higher-residual eigenpairs reveals numerous higher-order matches, forming lattice-like structures in the period--power plane. Residuals generally increase with power order, indicating that higher-order components are less robustly resolved by the finite-dimensional Koopman approximation. 
 } 
\label{Fig6}
\end{figure*}

Koopman theory provides a natural mechanism for structured spectral redundancy. If \(\phi_i\) and \(\phi_j\) are Koopman eigenfunctions with eigenvalues \(\lambda_i\) and \(\lambda_j\), then their product \(\phi_i\phi_j\), when well defined, is also an eigenfunction with eigenvalue \(\lambda_i\lambda_j\). Similarly, integer powers of a single eigenfunction generate eigenvalues \(\lambda_i^m\). Thus, once a fundamental oscillatory component is learned, k-EDMD may also produce eigenpairs corresponding to its harmonics or nonlinear products with other components. These eigenpairs need not represent independent climate modes; they may instead reflect algebraic structure generated by the Koopman eigenfunctions themselves. This viewpoint is closely related to the combination-mode interpretation used in previous Koopman analyses of ENSO, where interactions between ENSO and the annual cycle generate additional spectral components \citep{Froyland2021}. Here, we apply the same principle to the deseasonalized tropical Pacific anomaly spectrum, asking whether the dense k-EDMD cloud is organized by powers of robust residual-minimum components.

To test this organization, we use the matched k-EDMD modes as seeds and search the full spectrum for eigenpairs consistent with their integer powers. An eigenpair is assigned to the hierarchy of a seed mode when two criteria are met: the period implied by its eigenvalue lies within 10\% of the period expected for the corresponding power, and its eigenfunction subspace has correlation greater than 0.7 with the subspace generated by the corresponding eigenfunction power.

Figure~\ref{Fig6} visualizes the resulting power hierarchies. For each seed, we plot the periods of eigenpairs matched to its integer powers as a function of power order. The mode numbers in the legend indicate the rank of each matched k-EDMD seed after ordering the full spectrum by increasing residual. All seed modes associated with pseudospectral fundamental modes have low residual ranks: the largest rank is 126, corresponding to about 15.5\% of the full spectrum, while several key modes, including those with periods of 12.53, 5.75 and 4.73 years, lie within the first 20 modes. Thus, when discrete spectral counterparts to the pseudospectral modes are present, they are among the better-resolved k-EDMD eigenpairs.

Restricting the search to the lowest-residual 25\% or 50\% of the spectrum retains only a small set of low-order powers (Fig.~\ref{Fig6}, left). These low-residual matches define the best-resolved part of the hierarchy. Several coincide with pseudospectral residual minima (circled red markers), suggesting that some apparently fundamental frequencies are connected through low-order powers of slower Koopman components. Extending the search to the full spectral cloud reveals many additional higher-order matches (Fig.~\ref{Fig6}, right), forming lattice-like structures in the period--power plane. Although these higher-order components are generally less well resolved, their organization shows that the dense k-EDMD spectrum is not simply an unstructured collection of unrelated eigenpairs.

Residuals provide an operator-level ordering of this hierarchy. Low-order matches typically have smaller residuals and are therefore better resolved by the finite-dimensional Koopman approximation. As power order increases, residuals generally grow, indicating that higher-order components are less reliable and more sensitive to sampling variability and numerical approximation. The residual ordering therefore converts the dense spectrum into a graded hierarchy, ranging from robust fundamental modes and low-order dynamical relationships to weakly resolved high-residual structures.

These results recast spectral redundancy in k-EDMD as a form of Koopman-algebraic organization. Rather than representing independent climate degrees of freedom, a substantial part of the dense spectral cloud is arranged around integer powers, and potentially nonlinear combinations, of a compact set of fundamental Koopman components. Pseudospectral analysis identifies the robust backbone, while the full k-EDMD spectrum contains secondary components dynamically organized around it.

\subsubsection*{The Koopman backbone improves extended-lead Niño3.4 forecasts}

A compact Koopman backbone is useful only if it captures dynamical evolution, not merely spatial coherence or in-sample variance. We therefore test whether reduced-order models built from the pseudospectral fundamental modes retain out-of-sample forecast skill for tropical Pacific SST variability.

\begin{figure*}
 \centering
 \includegraphics[width=\textwidth]{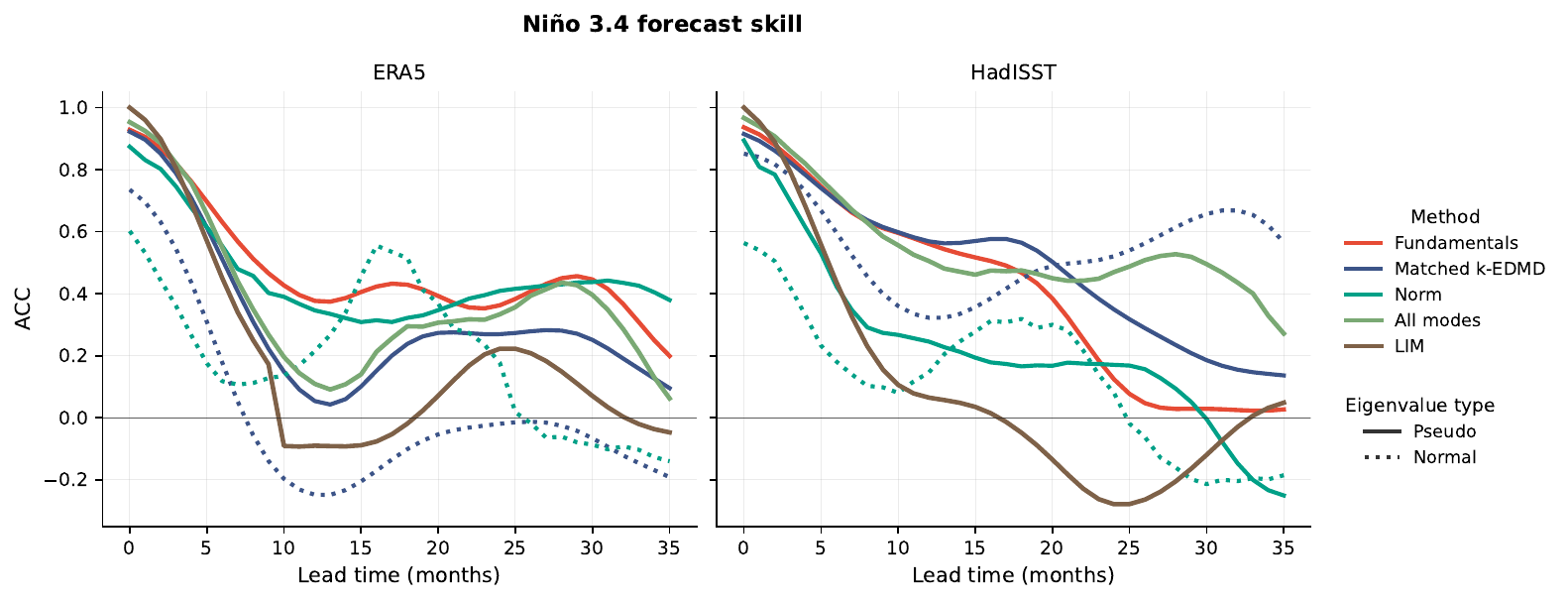}
 \caption{
\textbf{Out-of-sample Niño3.4 forecast skill from different modal representations.}
Anomaly correlation coefficient (ACC) of Niño3.4 forecasts as a function of lead time for ERA5 (left) and HadISST (right). Forecasts are generated from reduced-order models based on the pseudospectral fundamental modes, their matched k-EDMD counterparts, the full k-EDMD spectrum, modes selected by singular-value norm, and a Linear Inverse Model (LIM). Where applicable, solid lines denote forecasts using pseudoeigenfunctions, while dotted lines denote forecasts using standard k-EDMD eigenfunctions. Models based on the pseudospectral fundamental modes retain the strongest skill across the ENSO-relevant forecast range, particularly at intermediate lead times.
 }
 \label{Fig7}
\end{figure*}

Figure~\ref{Fig7} compares Niño3.4 forecast skill across modal representations. In both ERA5 and HadISST, reduced-order models built from the pseudospectral fundamental modes retain strong out-of-sample skill at extended lead times, with their clearest advantage at leads of approximately 8--18 months. This range is especially relevant for ENSO prediction and provides a stringent test of whether the selected modes encode dynamically persistent information rather than merely reconstructing training-period variability.

Forecasts based on the matched k-EDMD modes depend strongly on the eigenfunction representation. With standard k-EDMD eigenfunctions, skill decays rapidly with lead time, consistent with the sampling sensitivity of individual eigenpairs. Pseudoeigenfunctions substantially improve these forecasts, but the improvement is less reproducible across data products and training windows than for the pseudospectral fundamental modes. In HadISST, matched modes with pseudoeigenfunctions perform comparably to the fundamental modes up to lead times of approximately 20 months; in ERA5, their skill remains substantially lower. Thus, the frequencies and spatial patterns of matched modes can be dynamically meaningful, but their predictive value depends on how well the discrete k-EDMD spectrum resolves them in a given training realization. By contrast, pseudospectral fundamental modes provide a more stable forecasting basis because they are defined by robust residual minima rather than by individual eigenpairs in the dense spectral cloud.

Using the full k-EDMD spectrum does not improve forecast skill. Although the full spectral expansion can better match the initial state, it also retains many weakly resolved, high-residual components that degrade forecasts at intermediate and longer lead times, consistent with Lorenzo-Sánchez et al. \citep{Lorenzo2025}. Similarly, modes selected only by singular-value norm underperform the pseudospectral fundamental modes despite their large amplitudes. These comparisons show that amplitude-based selection is not sufficient: predictive skill requires modes that are both dynamically persistent and robustly resolved by the Koopman approximation.

The reduced-order Koopman forecasts also show higher skill than the Linear Inverse Model at intermediate lead times in both data products. This suggests that the nonlinear operator-theoretic representation captures predictive structure in tropical Pacific variability that is not fully represented by a linear stochastic model. The improvement is not simply a matter of model size, because forecasts using the full k-EDMD spectrum do not show the same advantage. Rather, the gain reflects the dynamical relevance and residual robustness of the pseudospectrally selected fundamental modes.

These results show that the Koopman backbone revealed by pseudospectral analysis is not merely an interpretable reconstruction basis. It captures dynamically persistent structures that support skillful ENSO forecasts outside the training period. By isolating robust residual-minimum modes from the dense and redundant k-EDMD spectrum, the framework provides a compact representation that is interpretable, dynamically organized and predictive.

\subsection*{Discussion}

Data-driven Koopman methods promise physically interpretable spectral representations of nonlinear climate variability, but in high-dimensional geophysical applications they often produce spectra that are too dense and sampling-sensitive to interpret mode by mode. Here we show that this apparent complexity has a reproducible operator-level organization. By combining k-EDMD with residual minimization and pseudospectral analysis, we identify a compact Koopman backbone of tropical Pacific SST variability: a set of robust residual-minimum modes that are stable across data products and sampling realizations, coherent in physical space, and predictive at ENSO-relevant lead times.

The main advance is to shift Koopman analysis from eigenpair selection to spectral organization. Individual k-EDMD eigenvalues and eigenvalue-density peaks vary under changes in data product or training window, whereas pseudospectral residual minima persist. These minima provide a dynamical consistency criterion for identifying well-resolved spectral components, rather than relying on amplitude, growth rate, variance contribution, or visual inspection of the eigenvalue cloud. This distinction is important for climate applications, where finite observational records, weak damping, multiscale variability, and near-continuous spectral structure can make individual learned eigenpairs unstable even when the underlying dynamical signal is reproducible.

The resulting backbone links three properties that are often optimized separately in data-driven climate modelling: robustness, interpretability, and prediction. The selected modes have coherent spatial structures spanning low-frequency tropical Pacific modulation, ENSO-band variability, central-Pacific ENSO-like patterns, and shorter quasi-biennial components. Their organization is consistent with previous operator-theoretic studies of Indo-Pacific and ENSO variability, while the present analysis adds an explicit residual and pseudospectral criterion for separating robust primary components from the surrounding finite-data spectral cloud. The same modes reconstruct most observed Niño3.4 variability and retain out-of-sample forecast skill, with strong advantages at ENSO-relevant lead times. Thus, the backbone is not only an interpretable decomposition of SST anomalies, but a reduced dynamical representation carrying predictive memory.

A second contribution is the interpretation of the remaining k-EDMD spectrum. Dense Koopman spectra are often treated as a practical obstacle: they contain many plausible modes, but only a small subset can be interpreted physically. Our results suggest a different view. Many additional eigenpairs are organized as integer powers, and potentially nonlinear combinations, of the residual-minimum components, consistent with the algebra of Koopman eigenfunctions. The dense spectrum therefore need not be interpreted as a collection of independent climate degrees of freedom. Instead, it contains structured redundancy around a compact set of robust components, with residuals providing an ordering from well-resolved fundamental modes to increasingly weakly resolved higher-order structures. This perspective explains why direct selection from the raw k-EDMD spectrum is unstable, while also showing that the spectral cloud contains useful dynamical information when viewed through residual and pseudospectral diagnostics.

These findings also clarify the physics-informed character of the framework. The imposed structure is not a prescribed climate equation or an externally specified mechanism, but the Koopman eigenvalue relation itself. Residual consistency acts as an operator-level physical constraint on the learned representation: components are retained when they are stable under perturbations, coherent in physical space, and consistent with the dynamics approximated by the learned operator. In this sense, residual and pseudospectral Koopman analysis provides a route from high-dimensional data-driven learning to reduced-order climate representations that are interpretable, dynamically organized, and predictive. This is particularly relevant for physics-informed machine learning, where a central challenge is to move beyond black-box accuracy toward learned structures that can be inspected, validated, and related to dynamical mechanisms.

Several limitations remain. The analysis focuses on tropical Pacific SST anomalies, where ENSO provides a strong and spatially coherent signal. It remains to be tested whether similar pseudospectral organization emerges in more strongly nonlinear, multivariate, or less directly ENSO-dominated fields, including subsurface temperature, wind stress, sea level, precipitation, and coupled ocean--atmosphere variables. The interpretation of decadal components also requires caution, because the observational record contains only a limited number of low-frequency cycles. Longer climate-model integrations, large ensembles, and controlled experiments will be needed to determine whether the identified backbone reflects a robust dynamical organization of the tropical Pacific system or a structure specific to SST over the available observational period.

A further open question concerns the physical meaning of the Koopman hierarchy. In a Fourier decomposition, harmonic structure is generated by fixed sinusoidal basis functions. In contrast, Koopman spectra may contain multiple families of powers and nonlinear products associated with state-dependent eigenfunctions. The number, spacing, residual ordering, and interactions of these families could provide a richer fingerprint of nonlinear dynamical organization than a single frequency cascade. Developing quantitative measures of this hierarchy may support mode selection, regime comparison, reconstruction, prediction, and the attribution of learned spectral components to specific ocean--atmosphere processes.

Overall, the results show that residual and pseudospectral Koopman analysis can transform dense and ambiguous learned climate spectra into compact, interpretable, and predictive dynamical representations. For tropical Pacific variability, this reveals a Koopman backbone that organizes ENSO-relevant timescales, reconstructs Niño3.4 variability, and improves extended-lead prediction. More broadly, the framework provides a physics-informed strategy for identifying robust low-dimensional structure in complex Earth-system data, especially in settings where interacting timescales and spectral redundancy obscure the dynamical content of learned representations.

\subsection*{Methods}

\subsubsection*{Datasets and data pre-processing}

To assess robustness across data products, we use monthly tropical Pacific sea-surface-temperature fields from ERA5 and HadISST. ERA5 is the ECMWF/Copernicus fifth-generation atmospheric reanalysis, available from 1940 onward and used here at \(0.25^{\circ} \times 0.25^{\circ}\) resolution \citep{hersbach2020era5}. HadISST provides globally complete monthly SST and sea-ice fields from 1870 onward on a \(1^{\circ} \times 1^{\circ}\) grid, based on historical in situ and satellite-era observations combined through reconstruction and interpolation procedures \citep{Rayner2003}. We use the common 1940--2025 period, with 1940--2010 used for training and the remaining years reserved for out-of-sample evaluation.

The analysis is restricted to the tropical Pacific domain 
\((31^\circ\mathrm{S}-32^\circ\mathrm{N},\,130^\circ\mathrm{E}-70^\circ\mathrm{W})\). 
After applying a land--sea mask, this domain contains 9,121 ocean grid points in HadISST and 146,586 in ERA5. At each grid point, we remove a linear trend to reduce the influence of long-term forced variability, apply a 3-month running mean to suppress high-frequency fluctuations, and form monthly anomalies by subtracting the 1940--2025 monthly climatology. The anomalies are then standardized by the local temporal standard deviation, so that the Koopman analysis is performed on basin-wide SST variability expressed in comparable units.

Before constructing the Koopman approximation, each dataset is projected onto its empirical orthogonal function (EOF) basis. Retaining the numerical rank of the anomaly matrix makes this projection equivalent to the full grid-point representation up to numerical precision, while substantially reducing computational cost \citep{navarra_2021}. The retained ranks are 189 for ERA5 and 746 for HadISST.

For the 1940--2025 analysis period, the 3-month running mean gives \(1018 (t)\) monthly states.  To account for memory effects in the reduced EOF-coordinate dynamics, we augment the state vector with 24 monthly time delays, corresponding to a maximum lag of two years. This gives \(t-24=994\) effective delayed states. Of these, the first 814 states \((m+1)\) are used for training the k-EDMD model and for the subsequent residual and pseudospectral analyses, while the remaining 180 states, corresponding to the final 15 years, are held out for forecast validation. The resulting data matrix is written as
\begin{equation*}\setlength\abovedisplayskip{6pt}\setlength\belowdisplayskip{6pt}\setlength\abovedisplayskip{6pt}\setlength\belowdisplayskip{6pt}
X = [\mathbf{x}_0,\mathbf{x}_1,\ldots,\mathbf{x}_m] \in \mathbb{R}^{r \times (m+1)},
\end{equation*}
where \(\mathbf{x}_i\) contains the retained EOF coefficients for the \(i\)-th monthly SST anomaly and \(r\) is the retained EOF rank. This EOF-coordinate representation restricted to 814 delayed training states (\(m = 813\)), forms the input to k-EDMD and to the subsequent residual and pseudospectral analyses.

\subsubsection*{Methodology}

We identify the fundamental Koopman components of ENSO variability by combining kernel-based Koopman approximation, residual diagnostics, pseudospectral analysis, and mode reconstruction. The aim is to extract physically interpretable, dynamically robust modes from the dense k-EDMD spectrum, using operator-level spectral resolution and stability as selection criteria rather than heuristic ranking of individual eigenpairs.

\vspace{1em}
\noindent\textbf{Koopman operator framework.}  
Let \( \mathbf{x}_{t+1} = \mathbf{F}(\mathbf{x}_t) \) denote a discrete dynamical system acting on the state space \( \mathcal{M} \subset \mathbb{R}^d \). The associated Koopman operator \( \mathcal{K} \) acts linearly on observables \( g : \mathcal{M} \to \mathbb{C} \) as
\begin{equation}\setlength\abovedisplayskip{6pt}\setlength\belowdisplayskip{6pt}\setlength\abovedisplayskip{6pt}\setlength\belowdisplayskip{6pt}
(\mathcal{K} g)(\mathbf{x}) = g(\mathbf{F}(\mathbf{x})).
\end{equation}
Despite the nonlinearity of \( \mathbf{F} \), \( \mathcal{K} \) is a linear (possibly infinite-dimensional) operator. Its spectral decomposition,
\begin{equation}\setlength\abovedisplayskip{6pt}\setlength\belowdisplayskip{6pt}\setlength\abovedisplayskip{6pt}\setlength\belowdisplayskip{6pt}
\mathcal{K}\phi_j = \lambda_j \phi_j, \qquad g(\mathbf{x}) = \sum_j \phi_j(\mathbf{x}) v_j,
\end{equation}
provides a representation of the dynamics in terms of eigenvalues \( \lambda_j \) (encoding decay and frequency) and Koopman modes \( v_j \) (spatial patterns associated with eigenfunction \( \phi_j \)).

\vspace{1em}
\noindent\textbf{Kernel-based Koopman approximation.}  
To approximate \( \mathcal{K} \) from data, we employ the kernel Extended Dynamic Mode Decomposition (k-EDMD) algorithm. Given snapshot matrices \( \mathbf{X} = [\mathbf{x}_0, \ldots, \mathbf{x}_{m}] \) and \( \mathbf{Y} = [\mathbf{y}_0, \ldots, \mathbf{y}_{m}] \), we construct Gram matrices in the associated reproducing kernel Hilbert space (RKHS):
\begin{equation}\setlength\abovedisplayskip{6pt}\setlength\belowdisplayskip{6pt}\setlength\abovedisplayskip{6pt}\setlength\belowdisplayskip{6pt}
\mathbf{G}_{xx} = k(\mathbf{X}, \mathbf{X}), \quad
\mathbf{G}_{yx} = k(\mathbf{Y}, \mathbf{X}), \quad
\mathbf{G}_{yy} = k(\mathbf{Y}, \mathbf{Y}),
\end{equation}
where \( k \) is a positive-definite kernel, in our case a Gaussian:
\begin{equation}\setlength\abovedisplayskip{6pt}\setlength\belowdisplayskip{6pt}\setlength\abovedisplayskip{6pt}\setlength\belowdisplayskip{6pt}
\label{gaussian_kernel}
\psi_i(\mathbf{x}) = \frac{1}{(\sqrt{2\pi}s)^{n/2}}\exp\left(-\frac{||\mathbf{x}-\mathbf{x}_i||^2}{2 s^2}\right)
\end{equation}
The Koopman operator is then approximated as
\begin{equation}\setlength\abovedisplayskip{6pt}\setlength\belowdisplayskip{6pt}\setlength\abovedisplayskip{6pt}\setlength\belowdisplayskip{6pt}
K = \mathbf{G}_{yx}\,\mathbf{G}_{xx}^{-1},
\end{equation}
acting in the finite subspace spanned by the kernel evaluations. Its eigendecomposition yields complex eigenvalues \( \lambda_j \), approximate eigenfunctions \( \Phi_j \) (evaluated at the data points), and associated Koopman modes \( V_j \).

\vspace{1em}
\noindent\textbf{Residual evaluation (kernelized version of ResDMD).}  
The accuracy of a finite-dimensional approximation can be quantified through residuals derived from the Koopman eigenvalue equation,
\begin{equation}\setlength\abovedisplayskip{6pt}\setlength\belowdisplayskip{6pt}\setlength\abovedisplayskip{6pt}\setlength\belowdisplayskip{6pt}
(\mathcal{K} - \lambda I)\phi = 0.
\end{equation}
The Residual DMD (ResDMD) algorithm measures the deviation of each approximate eigenpair \( (\lambda, \phi) \) from this equation.  
In the kernel setting, the normalized residual can be expressed as~\cite{Colbrook2024}:
\begin{equation}\setlength\abovedisplayskip{6pt}\setlength\belowdisplayskip{6pt}\setlength\abovedisplayskip{6pt}\setlength\belowdisplayskip{6pt}
\mathrm{res}^2(\lambda, \phi) =
\left| \sqrt{ \Re \!\left(
\frac{\operatorname{diag}(\mathbf{W}_2^* \mathbf{G}_{yy} \mathbf{W}_2)}
     {\operatorname{diag}(\mathbf{W}_2^* \mathbf{W}_2)}
\right) - |\lambda|^2 } \right|,
\label{residual_kernel}
\end{equation}
where \( \mathbf{W}_2 \) denotes the matrix of left eigenvectors of \( K \) in the coefficient space, and \( \mathbf{G}_{yy} \) approximates the action of \( K^* K \) in the RKHS.  
This metric quantifies how closely each computed eigenfunction satisfies the Koopman eigenvalue relation under the chosen kernel. Eigenpairs with low residuals correspond to spectrally robust, physically meaningful dynamics.

\vspace{1em}
\noindent\textbf{Pseudospectral analysis.}  
Because the Koopman operator is generally non-normal, small perturbations in the data or kernel approximation can lead to large changes in its spectrum. To evaluate spectral robustness, we compute the \textit{pseudospectrum} \( \Lambda_\varepsilon(\mathcal{K}) = \{ \mu \in \mathbb{C} : \|( \mathcal{K} - \mu I)^{-1} \| > \varepsilon^{-1} \} \) of the estimated Koopman operator.
In discrete form, for each complex number \( \mu \) on the unit circle, we consider the perturbed operator
\begin{equation}\setlength\abovedisplayskip{6pt}\setlength\belowdisplayskip{6pt}\setlength\abovedisplayskip{6pt}\setlength\belowdisplayskip{6pt}
M(\mu) = \mathbf{G}_{yy} - \mu \bar{\mathbf{A}} - \bar{\mu}\mathbf{A} + \|\mu^2\|\mathbf{G}_{xy},
\end{equation}
where \( \mathbf{A} \) is the auxiliary matrix defined as $(\mathbf{G}_{xx} + n \epsilon I )^{-1}\mathbf{G}_{y x}$.  
The smallest singular value \( \sigma_{\min}(M(\mu)) \) indicates how close \( \mu \) is to being an eigenvalue. We normalize by a scaling matrix \( \mathbf{S_Q} \) and define the pseudospectral residual as
\begin{equation}\setlength\abovedisplayskip{6pt}\setlength\belowdisplayskip{6pt}\setlength\abovedisplayskip{6pt}\setlength\belowdisplayskip{6pt}
\mathrm{Res}(\mu) = \sqrt{\sigma_{\min}(\mathbf{S_Q} \, M(\mu) \, \mathbf{S_Q})}.
\end{equation}
Evaluating \( \mathrm{Res}(\mu) \) along the unit circle yields a residual spectrum in the frequency domain, where minima identify the \textit{fundamental frequencies} most strongly supported by the dynamics.

\vspace{1em}
\noindent\textbf{Mode reconstruction via spectral filtering.}  
For each selected fundamental frequency \( \mu_t \), we reconstruct the corresponding Koopman mode directly from the time series of SST fields \( \mathbf{x}_0, \ldots, \mathbf{x}_{m+1} \) without relying on the kernel dictionary. The mode \( \boldsymbol{\rho}_u \) is obtained via time-domain averaging:
\begin{equation}\setlength\abovedisplayskip{6pt}\setlength\belowdisplayskip{6pt}\setlength\abovedisplayskip{6pt}\setlength\belowdisplayskip{6pt}
\boldsymbol{\rho}_u = \mathrm{Re}\left\{ \frac{1}{N} \sum_{k=0}^{m+1} \mu_t^{-k}\,\mathbf{x}_k \right\}.
\label{summation_modes}
\end{equation}
This \textit{Koopman spectral filtering} isolates coherent spatial structures oscillating with frequency and decay rate determined by \( \mu_t \), providing a complementary, nonparametric reconstruction of the mode \cite{Colbrook2026KoopmanLearning, Mohr2014Construction}.

\vspace{1em}
\noindent\textbf{Identification of matched components}  
We first compare the reconstructed fundamental modes with the dictionary-based k-EDMD modes. For each reconstructed mode \( \boldsymbol{\rho}_u \) associated with a fundamental frequency \( \mu_t \), we identify the closest k-EDMD eigenvalue \( \lambda_k \) in the complex plane and compute the spatial anomaly correlation coefficient (ACC) between their corresponding spatial patterns. A threshold of 0.7 is adopted to define \textit{matched} modes, which are retained as the k-EDMD equivalents to the fundamental components of the system.  
This matching step links the dynamically most coherent frequencies obtained from the pseudospectral analysis with the data-driven Koopman modes recovered by k-EDMD.

\vspace{0.5em}
To further exploit the information contained in the fundamental frequencies, we also compute the corresponding \textit{pseudo-eigenfunctions}. These functions are obtained by projecting the observable time series onto the complex exponential basis associated with each target eigenvalue, thereby approximating the Koopman eigenfunction on the data manifold. For a given observable sequence \( g_t = g(\mathbf{x}_t) \) and target eigenvalue \( \mu_t \), the pseudo-eigenfunction \( \tilde{\phi}_{\mu_t} \) is defined as
\begin{equation}\setlength\abovedisplayskip{6pt}\setlength\belowdisplayskip{6pt}\setlength\abovedisplayskip{6pt}\setlength\belowdisplayskip{6pt}
\tilde{\phi}_{\mu_t}(\mathbf{x}_t) = \frac{1}{N} \sum_{k=0}^{m+1} \mu_t^{-k} \, g(\mathbf{x}_{t+k}),
\label{pseudo_eig}
\end{equation}
which can be interpreted as a discrete projection of the observable onto the spectral subspace associated with \( \mu_t \). Unlike the kernel-based eigenfunctions \( \Phi_j \) from k-EDMD, these pseudo-eigenfunctions are obtained directly from the data and depend solely on the target frequencies identified in the residual spectrum.

\vspace{0.5em}
Both the k-EDMD eigenfunctions and the pseudo-eigenfunctions are then used for data reconstruction and forecasting. The k-EDMD modes provide a dictionary-based decomposition in the RKHS, while the pseudo-eigenfunctions yield a complementary, nonparametric representation derived from direct temporal filtering. Their joint use allows for a more complete characterization of the underlying dynamics and an evaluation of the predictive skill associated with each target frequency.  

This full procedure is applied independently to ERA5 and HadISST datasets. Common fundamental components and pseudo-eigenfunctions across both datasets are interpreted as physically robust features of the ENSO dynamics.

\vspace{1em}
\noindent\textbf{Identification of nonlinear mode hierarchies.}  
Finally, we investigate the algebraic structure of the Koopman spectrum. Because the Koopman operator is closed under composition, products and powers of eigenfunctions are themselves eigenfunctions with eigenvalues \( \lambda^n \). We therefore search for eigenvalue pairs satisfying
$
\lambda_k \approx \lambda_f^n
$
for integer \( n > 1 \), within specified period tolerances.  
If, in addition, the subspace angle between their eigenfunctions is small (or their correlation high), the pair is accepted as a power relation. A 10\% error around the period and a 0.7 ACC threshold for the eigenfunctions have been used in this study. This analysis tests whether the observed Koopman spectrum exhibits a lattice-like hierarchy, potentially revealing a small set of fundamental modes generating higher-order harmonics of ENSO variability.

\backmatter

\bmhead{Supplementary information}

A complementary PDF file containing the supplementary information, including additional figures and tables, is provided with this article.

\bmhead{Data and code availability} The data and code supporting the findings of this study are available in the dedicated GitHub repository: https://github.com/Paulalo95/ENSO-Fundamental-Modes.

\bmhead{Contributions} P.L.S. designed and performed the numerical experiments, carried out the data analysis, interpreted the results, and wrote the initial draft. M.J.C. helped shape the study’s mathematical and methodological direction, with particular contributions to the residual and pseudospectral Koopman framework and the identification of key computations. M.J.C. also contributed to the writing and revision of the manuscript. A.N. provided climate-science supervision and expertise on tropical Pacific variability and ENSO, and helped interpret the physical significance of the results. All authors contributed to discussions, commented on earlier versions of the manuscript, and approved the final manuscript.

\renewcommand{\baselinestretch}{0.7}
\bibliographystyle{abbrv}
\footnotesize
\bibliography{sn-bibliography}

\clearpage
\onecolumn

\appendix
\section*{Supplementary Material}
\addcontentsline{toc}{section}{Supplementary Material}

\noindent
This document contains supplementary figures, tables, and methodological details supporting the main manuscript.\\

The supplementary material provides additional evidence for the residual and pseudospectral organization of the dense k-EDMD spectrum. Table~S1 summarizes the 19 pseudospectral fundamental periods identified in ERA5, their HadISST counterparts, and the matched k-EDMD eigenmodes where present, together with spatial correlations and assigned dynamical families. The strong agreement across data products for most modes shows that the selected frequencies and associated spatial structures are reproducible, while the matched k-EDMD modes demonstrate that these pseudospectral components are embedded within the original dense Koopman spectrum.

Figure~S1 examines the sensitivity of the k-EDMD spectrum to small changes in the training window. Although the detailed eigenvalue clouds and eigenvalue-density peaks vary substantially across adjacent 65-year windows, the pseudospectral residual curves remain highly consistent, with local minima occurring at nearly the same periods. This supports the use of residual minima, rather than raw eigenvalue-density peaks, as robust indicators of dynamically meaningful frequencies.

Figures~S2 and~S3 show the full set of ERA5 pseudospectral fundamental modes. The modes display coherent tropical Pacific SST anomaly structures across periods ranging from approximately 1.5 years to decadal timescales. Longer-period modes exhibit broader basin-scale and off-equatorial patterns, ENSO-band modes are concentrated along the central and eastern equatorial Pacific cold tongue, and shorter-period modes become increasingly localized. Together, these supplementary results support the interpretation of the selected modes as a compact, physically interpretable Koopman backbone of tropical Pacific variability.

\section{Supplementary Tables}
\label{sec:supp_tables}

\begin{table}[h]
\centering
\small
\caption{
\textbf{Summary of pseudospectral fundamental modes.}
ERA5 fundamental periods (identified from pseudospectral residual minima) are present in the first column. Columns 2-3 report the periods of the equivalent HadISST pseudospectral modes and matched k-EDMD modes, where identified, together with theit respective spatial correlations (in parentheses). The fifth column lists the assigned physical or dynamical family of each mode.
}
\label{tab:fundamental_modes}
\begin{tabular}{ccccc}
\hline
\textbf{Mode} & \textbf{ERA5 period (yr)} & \textbf{HadISST} & \textbf{k-EDMD} & \textbf{Dynamical family} \\
\hline
1  & 12.63 & 12.44 (0.88) & 12.53 (0.91) & Low-frequency Pacific modulation \\
2  & 8.77  & 9.59 (0.81) & 8.98 (0.85) & Low-frequency Pacific modulation \\
3  & 6.56  & 6.56 (0.77) & -- & ENSO-band variability \\
4  & 5.67  & 5.67 (0.92) & 5.75 (0.93) & ENSO-band variability \\
5  & 4.96  & 4.96 (0.93) & 4.73 (0.70) & ENSO-band variability \\
6  & 4.19  & 4.19 (0.91) & 3.99 (0.76) & ENSO-band variability \\
7  & 3.62  & 3.62 (0.97) & 3.72 (0.94) & ENSO-band variability \\
8  & 3.03  & -- & -- & ENSO-band variability \\
9  & 2.87  & 2.93 (0.73) & 2.88 (0.92) & Power 2 of Mode 4 \\
10 & 2.56  & -- & 2.62 (0.80) & Fast ENSO variability \\
11 & 2.44  & 2.44 (0.85) & 2.43 (0.84) & Fast ENSO variability \\
12 & 2.31  & 2.31 (0.77) & 2.33 (0.79) & Power 2 of Mode 5 \\
13 & 2.14  & 2.16 (0.71) & 2.17 (0.88) & High frequency modes \\
14 & 1.97  & 1.97 (0.86) & 1.95 (0.81) & Power 3 of Mode 4\\
15 & 1.90  & 1.90 (0.82) & 1.93 (0.83) & Power 3 of Mode 4\\
16 & 1.77  & -- & 1.75 (0.88) & High frequency modes \\
17 & 1.66  & 1.66 (0.73)  & -- & High frequency modes \\
18 & 1.54  & -- & -- & High frequency modes \\
19 & 1.48  & 1.48 (0.91) & 1.49 (0.80) & High frequency modes \\
\hline
\end{tabular}
\end{table}

\newpage
\section{Supplementary Figures}
\label{sec:supp_figures}

\begin{figure}[h]
 \centering
 \includegraphics[width=0.95\textwidth]{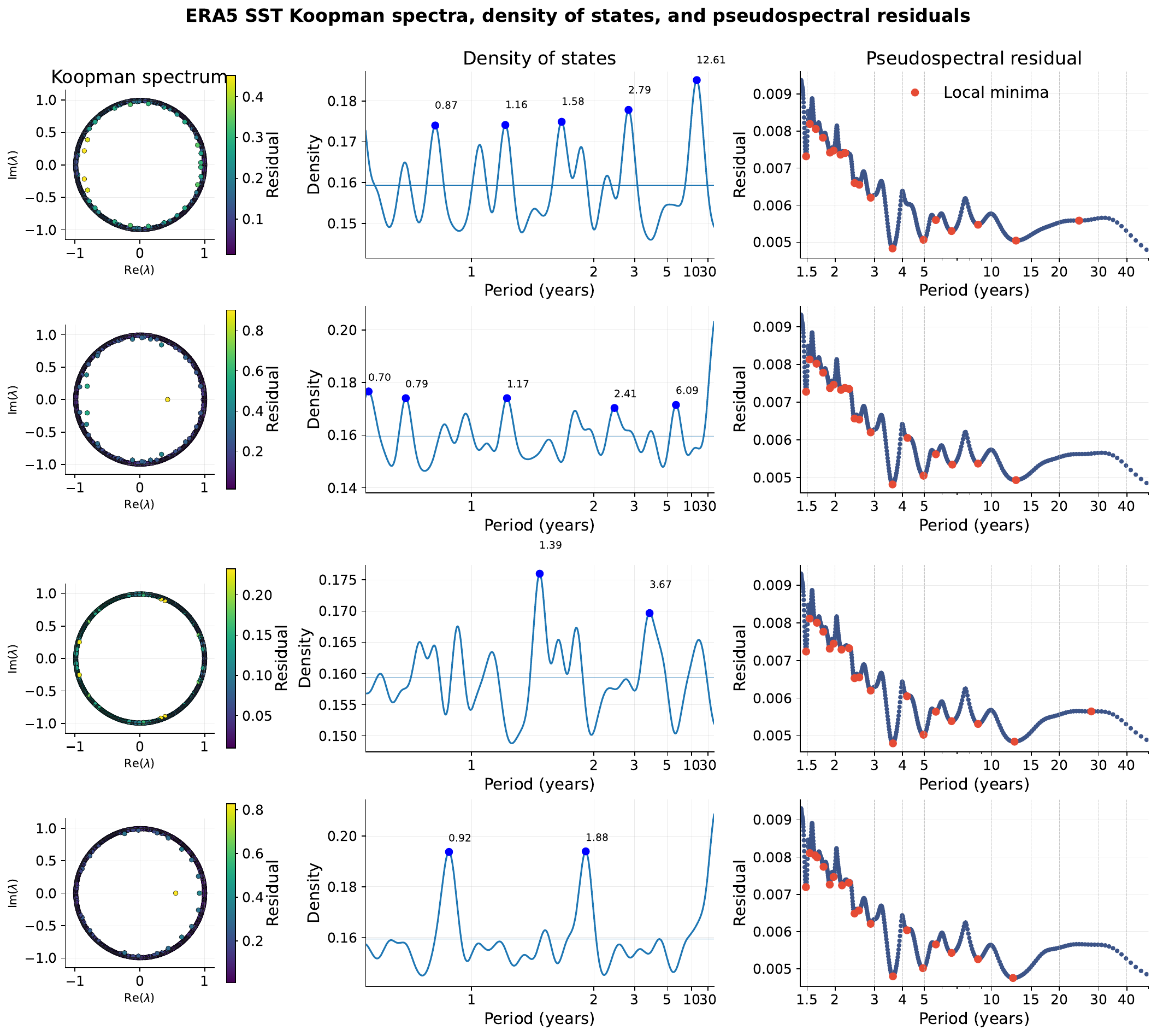}
 \caption{\textbf{Sensitivity of k-EDMD spectra and pseudospectral residual across nearby training windows.}
k-EDMD spectra, eigenvalue-density estimates, and pseudospectral residual curves for 65-year ERA5 tropical Pacific SST training windows starting in 1940, 1941, 1942, and 1943 (from top to bottom). In each row, the left panel shows k-EDMD eigenvalues coloured by residual norm, the middle panel shows the eigenvalue-density estimate as a function of period, and the right panel shows the unit-circle pseudospectral residual, with local minima marked in red. Although the eigenvalue clouds and density peaks vary across adjacent windows, the residual minima remain nearly unchanged, indicating that the pseudospectral residual structure is more stable than individual eigenpairs or raw density peaks.
} 
\label{AFig1}
\end{figure}

\begin{landscape}
\begin{figure}[p]
    \centering
    \includegraphics[width=0.85\linewidth]{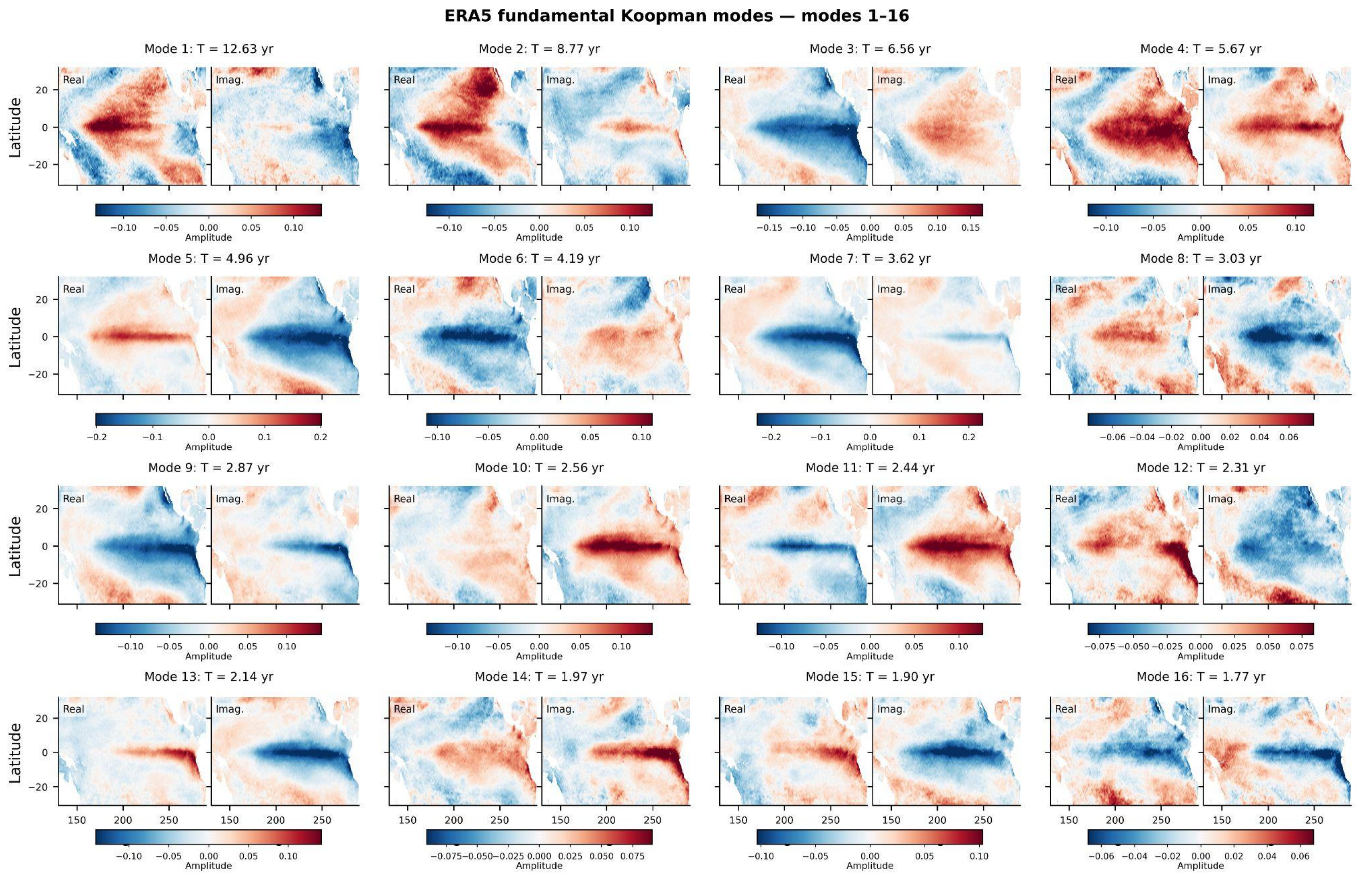}
    \caption{
\textbf{Pseudospectral fundamental Koopman modes in ERA5.}
Real and imaginary components of the fundamental Koopman modes identified from the pseudospectral analysis of ERA5 tropical Pacific SST variability. 
Each mode is shown as a pair of spatial patterns, with the real component on the left and the imaginary component on the right. 
The oscillation period $T$ associated with each mode is reported above the corresponding pair.  
Red and blue shading indicate opposite signs of the modal amplitude. 
The modes display coherent tropical Pacific structures across a range of interannual periods, highlighting the spatial organization of the fundamental components selected by the pseudospectral criterion.
}
    \label{fig:supp_landscape_1}
\end{figure}
\end{landscape}

\clearpage

\begin{landscape}
\begin{figure}[p]
    \centering
    \includegraphics[width=0.85\linewidth]{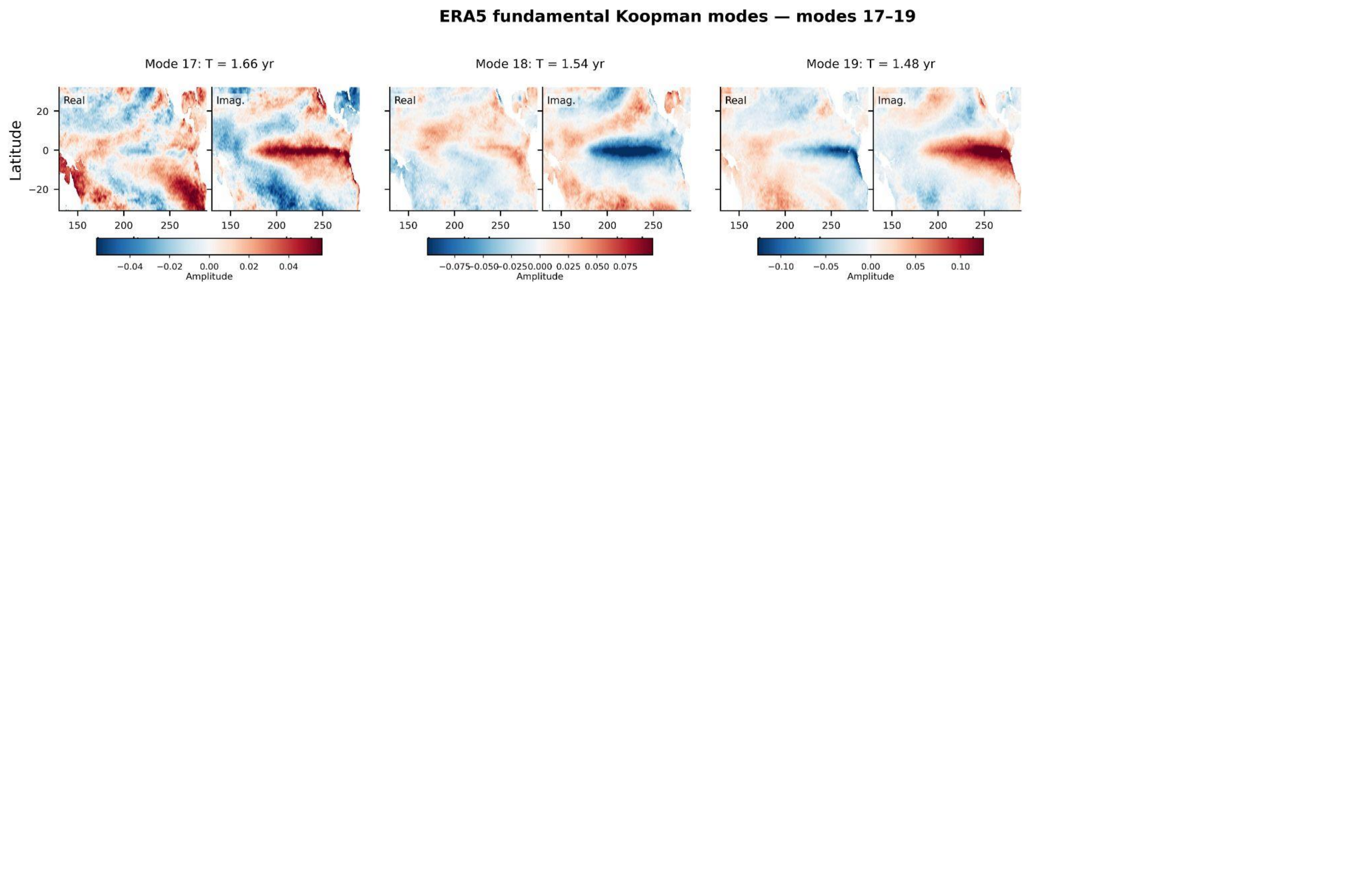}

    \caption{
\textbf{Pseudospectral fundamental Koopman modes in ERA5.}
Same as Figure S2 but for the Modes 17-19.
}
    
\label{fig:supp_fundamental_modes_ERA5_1_16}
    \label{fig:supp_landscape_2}
\end{figure}
\end{landscape}

\end{document}